\documentclass[12pt]{article}

\newbox\squ  % box character for ends of proofs
\setbox\squ=\hbox{\vrule width.3pt
            \vbox{\hrule height.3pt width.4em\kern1ex\hrule height.3pt}%
            \vrule width.3pt}
\def\endproof{%
  \ifmmode\eqno\copy\squ\medskip\else{\unskip\nobreak\hfil%
    \penalty50\hskip2em\hbox{}\nobreak\hfil\copy\squ
    \parfillskip=0pt \finalhyphendemerits=0\penalty-100\medskip}
  \fi} 

\usepackage{amsmath}
\usepackage{amssymb}

\begin{document}

\newcommand{\non}{\nonumber}
\newcommand{\wt}{\widetilde}
\newcommand{\wh}{\widehat}
\newcommand{\ot}{\otimes}
\newcommand{\g}{\mathfrak g}
\newcommand{\ts}{\,}
\newcommand{\U}{ {\rm U}}
\newcommand{\CC}{ {\rm C}}
\newcommand{\Z}{ {\rm Z}}
\newcommand{\R}{ {\rm R}}
\newcommand{\J}{ {\rm J}}
\newcommand{\Y}{ {\rm Y}}
\newcommand{\C}{\mathbb{C}}
\newcommand{\ZZ}{\mathbb{Z}}
\newcommand{\sdet}{ {\rm sdet}\ts}
\newcommand{\h}{\mathfrak h}
\newcommand{\gl}{\mathfrak{gl}}
\newcommand{\oa}{\mathfrak{o}}
\newcommand{\spa}{\mathfrak{sp}}
\newcommand{\Proof}{\noindent{\it Proof.}\ \ }   
\renewcommand{\theequation}{\arabic{section}.\arabic{equation}}  

\newtheorem{thm}{Theorem}[section]
\newtheorem{prop}[thm]{Proposition}
\newtheorem{cor}[thm]{Corollary}
\newtheorem{defin}[thm]{Definition}
\newtheorem{lem}[thm]{Lemma}

\newcommand{\bth}{\begin{thm}}
\renewcommand{\eth}{\end{thm}}
\newcommand{\bpr}{\begin{prop}}
\newcommand{\epr}{\end{prop}}
\newcommand{\bco}{\begin{cor}}
\newcommand{\eco}{\end{cor}}

\title{\Large\bf A weight basis for representations\\
of even orthogonal Lie algebras}
\author{{\sc A. I. Molev}\\[15pt]
School of Mathematics and Statistics\\
University of Sydney,
NSW 2006, Australia\\
alexm@maths.usyd.edu.au}%\\[30pt]
%Research Report 99--??}

\date{} % February 1999
\maketitle

\vspace{7 mm}

\begin{abstract}
A weight basis for each finite-dimensional irreducible representation
of the orthogonal Lie algebra $\oa(2n)$ is constructed.
The basis vectors are parametrized by the $D$-type 
Gelfand--Tsetlin patterns. The basis is consistent
with the chain of subalgebras $\g_1\subset\cdots \subset\g_n$,
where $\g_k=\oa(2k)$. Explicit formulas for the matrix elements 
of generators of $\oa(2n)$ in this basis are given. 
The construction is based on
the representation theory of the Yangians and extends
our previous results for the symplectic Lie algebras.

\end{abstract}

\newpage

\section{Introduction}
\setcounter{equation}{0}

In their pioneering works \cite{gt:fdu} and  \cite{gt:fdo}
Gelfand and Tsetlin proposed a combinatorial method
to explicitly construct representations of the classical Lie algebras.
For each finite-dimensional irreducible
representation of the
general linear Lie algebra $\gl(N)$ and 
the orthogonal Lie algebra $\oa(N)$
they gave a parametrization
of basis vectors and provided explicit
formulas for the matrix
elements of generators
of the Lie algebras in the basis. 
Derivations of the matrix element formulas 
in the orthogonal case are given in 
\cite{ph:lr,w:ro}; see also \cite{g:wc}.
A number of different approaches
to the problem of constructing representation
bases for simple Lie algebras has been developed;
see \cite{m:br} for more references. Note also recent results by
Donnelly~\cite{d:ec} and Littelmann~\cite{l:cc}.
In \cite{d:ec}
explicit combinatorial
constructions of the fundamental representations
of the $B$ and $C$ series Lie algebras 
and of their $q$-analogs are given;
in \cite{l:cc} monomial bases parametrized by patterns of
Gelfand--Tsetlin type are constructed for all simple
complex Lie algebras.
In \cite{m:br} an analog of
the Gelfand--Tsetlin basis for the
symplectic Lie algebras $\spa(2n)$ is constructed 
and explicit formulas
for the matrix elements of generators
of $\spa(2n)$ in this basis are given.
Bases for the finite-dimensional irreducible
representations of the classical Lie algebras
of $B$, $C$, and $D$ series can be constructed in
a uniform manner with the use of the representation 
theory of the Yangians, as in \cite{m:br}.
Here we extend the results of \cite{m:br}
to the case of the $D$ series and hope to treat the remaining
$B$ series case in a forthcoming publication.

Our basis for $\oa(2n)$
is different from that of Gelfand and Tsetlin \cite{gt:fdo}.
Their basis is consistent with the chain of subalgebras
\begin{equation}\label{chainGT}
\oa(2)\subset\oa(3)\subset\cdots\subset\oa(N).
\nonumber
\end{equation}
The reductions $\oa(k)\downarrow\oa(k-1)$ are multiplicity-free
which makes the basis orthogonal with respect to a
natural contravariant bilinear form. However, the basis vectors
are not weight vectors with respect to the
Cartan subalgebra of $\oa(N)$. 
To get a weight (although non-orthogonal) basis we consider 
the following chain instead:
\begin{equation}\label{chain}
\oa(2)\subset\oa(4)\subset\cdots\subset\oa(2n)
\nonumber
\end{equation}
so that all the subalgebras
belong to the $D$ series.

The reduction
$\oa(2n)\downarrow\oa(2n-2)$
is not multiplicity free.
This means that the subspace $V(\lambda)^{+}_{\mu}$
of $\oa(2n-2)$-highest vectors of a weight
$\mu$ in an $\oa(2n)$-module $V(\lambda)$ is not 
necessarily one-dimensional. However, this
space turns out to possess a natural structure
of an irreducible representation of a large associative algebra
$\Y^+(2)$ called the twisted Yangian (introduced by
Olshanski in \cite{o:ty}) 
and can also be equipped
with an action of the $\gl(2)$-Yangian $\Y(2)$. This allows us to
construct a Yangian Gelfand--Tsetlin basis in $V(\lambda)^{+}_{\mu}$
associated with an inclusion $\Y(1)\subset\Y(2)$; see 
\cite{m:gt,nt:yg,nt:ry}.

Our calculations are based on the relationship between
the twisted Yangian
$\Y^+(2)$ and the transvector
algebra $\Z(\g_n,\g_{n-1})$, $\g_n=\oa(2n)$. The transvector algebras
(they are sometimes called the Mickelsson algebras or $S$-algebras)
are studied in detail in \cite{z:gz,z:it}.

Although the constructions of the bases are very similar
for the orthogonal and symplectic cases, there is a slight
difference in the calculation of the
matrix elements of generators of the Lie algebras in the basis.
The Lie algebra $\spa(2n)$ contains the ``second diagonal"
generators $F_{-k,k}=2E_{-k,k}$ where the $E_{ij}$ 
denote the standard generators of
$\gl(2n)$ (see Section~\ref{sec:npr} below). 
The action of these elements in the basis
is rather simple and can be easily found. 
However, their counterparts
do not exist in the orthogonal
case. Instead, there are second degree elements $\Phi_{-k,k}$
of the universal enveloping algebra
$\U(\g_n)$ which belong
to the centralizer of $\g_{k-1}$ in $\U(\g_k)$ and play the role
similar to that of the elements $F_{-k,k}$
in the symplectic case.

\section{Notations and preliminary results}\label{sec:npr}
\setcounter{equation}{0}

We shall enumerate the rows and columns of $2n\times 2n$-matrices over
$\C$ by the indices $-n,\dots,-1,1,\dots,n$.
We let the $E_{ij}$, $i,j=-n,\dots,n$ 
denote the standard
basis of the Lie algebra $\gl(2n)$. 
We shall also assume
throughout the paper that the index $0$ is skipped in
a sum or in a product. Introduce the elements
\begin{equation}\label{Fij}
F_{ij}=E_{ij}-E_{-j,-i}. 
\end{equation}
We have $F_{-j,-i}=-F_{ij}$. In particular, $F_{-i,i}=0$ for all $i$.
The orthogonal Lie algebra $\g_n:=\oa(2n)$
can be identified
with the subalgebra in $\gl(2n)$
spanned by the elements $F_{ij}$, $i,j=-n,\dots,n$. 

The subalgebra $\g_{n-1}$ is spanned by the elements (\ref{Fij}) with the
indices $i,j$ running over the set $\{-n+1,\dots,n-1\}$. 
Denote by $\h=\h_n$ the diagonal Cartan subalgebra in $\g_n$. 
The elements $F_{11},\dots,F_{nn}$ form a basis of $\h$. 

The finite-dimensional irreducible representations of $\g_n$
are in a one-to-one correspondence with $n$-tuples
$\lambda=(\lambda_1,\dots,\lambda_n)$ 
where all the entries $\lambda_i$ are simultaneously
integers or half-integers 
(elements of the set $\frac12+\ZZ$)
and the following
inequalities hold:
\begin{equation}\label{inD}
-|\lambda_1|\geq\lambda_2\geq\cdots\geq \lambda_n.
\non
\end{equation}
Such an $n$-tuple $\lambda$ is called the highest weight
of the corresponding representation which
we shall denote by $V(\lambda)$.
It contains a unique, up to a multiple, nonzero vector $\xi$
(the highest vector) such that
$
F_{ii}\ts\xi=\lambda_i\ts\xi
$
for $i=1,\dots,n$ and
$
F_{ij}\ts\xi=0
$
for $-n\leq i<j\leq n$.

Denote by $V(\lambda)^+$ the subspace of $\g_{n-1}$-highest vectors
in $V(\lambda)$:
\begin{equation}\label{subhv}
V(\lambda)^+=\{\eta\in V(\lambda)\ |\ F_{ij}\ts \eta=0,
\qquad -n<i<j<n\}.
\non
\end{equation}
Given a $\g_{n-1}$-highest weight
$\mu=(\mu_1,\dots,\mu_{n-1})$ we denote by $V(\lambda)^+_{\mu}$
the corresponding weight subspace in $V(\lambda)^+$:
\begin{equation}\label{hvwmu}
V(\lambda)^+_{\mu}=\{\eta\in V(\lambda)^+\ |\ F_{ii}\ts\eta=
\mu_i\ts\eta,\qquad i=1,\dots,n-1\}.
\non
\end{equation}

Consider the extension 
of the universal enveloping algebra $\U(\g_n)$
\begin{equation}\label{ext}
\U'(\g_n)=\U(\g_n)\ot_{\U(\h)} \R(\h),
\non
\end{equation}
where $\R(\h)$ is the field of fractions of the commutative algebra
$\U(\h)$. Let $\J$
denote the left ideal in $\U'(\g_n)$ generated by 
the elements $F_{ij}$ with $-n<i<j<n$.
Set
\begin{equation}\label{tra}
\Z(\g_n,\g_{n-1})=\{x\in \U'(\g_n)/\J\ |\  F_{ij}\ts x\equiv0,
\quad -n<i<j<n\}.
\non
\end{equation}
Then $\Z(\g_n,\g_{n-1})$
is an algebra with the multiplication inherited from $\U'(\g_n)$.
We call it the {\it transvector algebra\/};
see \cite{z:gz,z:it} for further details. 

Set
\begin{equation}\label{fiD}
f_i=F_{ii}-i+1,\qquad f_{-i}=-f_i
\non
\end{equation}
for $i=1,\dots,n$.
Let $p$ denote the
{\it extremal projection\/}
for the Lie algebra $\g_{n-1}$; see \cite{ast:po,z:it}.
The projection $p$ naturally acts in the space
$\U'(\g_n)/\J$ and its image coincides with $\Z(\g_n,\g_{n-1})$.
The elements
\begin{equation}%\label{gener}
pF_{ia},\qquad\qquad a=-n,n,\qquad i=-n+1,\dots,n-1
\non
\end{equation}
are generators of $\Z(\g_n,\g_{n-1})$ \cite{z:it}. 
They can be given by
the following explicit formulas (modulo $\J$):
\begin{equation}%\label{pFia}
pF_{ia}=\sum_{i>i_1>\cdots>i_s>-n}
F_{ii_1}F_{i_1i_2}\cdots F_{i_{s-1}i_s}F_{i_sa}
\frac{1}{(f_i-f_{i_1})\cdots (f_i-f_{i_s})},
\non
%\\ \label{pFai} pF_{ai}&=\sum_{i<i_1<\cdots<i_s<n}
%F_{i_1i}F_{i_2i_1}\cdots F_{i_si_{s-1}}F_{ai_s}
%\frac{1}{(f_i-f_{i_1})\cdots (f_i-f_{i_s})},\non
\end{equation}
where $s=0,1,\dots$ (it is assumed that index $0$ is excluded
in the sum). We shall use the normalized
generators of $\Z(\g_n,\g_{n-1})$ defined by
\begin{align}\label{zia}
z_{ia}&=pF_{ia}(f_i-f_{i-1})\cdots\wh{(f_i-f_{-i})}\cdots(f_i-f_{-n+1}),\\
%\label{zai}
z_{ai}&=pF_{ai}(f_i-f_{i+1})\cdots\wh{(f_i-f_{-i})}\cdots(f_i-f_{n-1}),
\non
\end{align}
where the hats indicate the factors to be omitted if they occur.
We obviously have $z_{ai}=(-1)^{n-i}\ts z_{-i,-a}$.
The following equivalent formula holds for $z_{ai}$:
\begin{multline}\label{zaieq}
z_{ai}=(f_i-f_{i+1})\cdots\wh{(f_i-f_{-i})}\cdots(f_i-f_{n-1})\\
{}\times\sum_{n>i_1>\cdots>i_s>i}\frac{1}{(f_i-f_{i_1})\cdots (f_i-f_{i_s})}
F_{ai_1}F_{i_1i_2}\cdots F_{i_{s-1}i_s}F_{i_si}.
\end{multline}

The elements $z_{ia}$ and $z_{ai}$ 
naturally act in the space 
$V(\lambda)^+$ and are called
the {\it raising\/} and {\it lowering operators\/}.
One has for $i=1,\dots,n-1$:
\begin{equation}%\label{rlo}
z_{ia}:V(\lambda)^+_{\mu}\to V(\lambda)^+_{\mu+\delta_i},\qquad
z_{ai}:V(\lambda)^+_{\mu}\to V(\lambda)^+_{\mu-\delta_i},\non
\end{equation}
where $\mu\pm\delta_i$ is obtained from $\mu$ by replacing $\mu_i$
with $\mu_i\pm 1$.

Note the following relations between these operators; cf.  \cite{z:it}.
For $a,b\in\{-n,n\}$ and $i+j\ne 0$ one has
\begin{equation}%\label{relcom}
z_{aj}z_{bi}(f_i-f_j+1)=z_{bi}z_{aj}(f_i-f_j)
+z_{ai}z_{bj}.
\non
\end{equation}
In particular, $z_{ai}$ and $z_{aj}$ commute for $i+j\ne 0$.
One easily verifies that $z_{ai}$ and $z_{bi}$ also commute for
all $a,b$.
We shall use the following element which can be checked to belong
to the algebra $\Z(\g_n,\g_{n-1})$
\begin{equation}\label{zn-n}
z_{n,-n}=
\sum_{n>i_1>\cdots>i_s>-n}
F_{ni_1}F_{i_1i_2}\cdots F_{i_s,-n}\ts
\frac{(f_{n}-f_{j_1})\cdots (f_{n}-f_{j_k})}{2f_n},
\end{equation}
where $s=1,2,\dots$ and $\{j_1,\dots,j_k\}$ is the complement
to the subset $\{i_1,\dots,i_s\}$ in
$\{-n+1,\dots,n-1\}$. There is an equivalent formula for $z_{n,-n}$
which can either be proved directly (cf. \cite[Section~2]{m:br}),
or can be deduced from (\ref{eval}) below
(use the fact that $Z_{n,-n}(g_n)=Z_{n,-n}(-g_n)$):
\begin{equation}%\label{zn-neq}
z_{n,-n}=
\sum_{n>i_1>\cdots>i_s>-n}
F_{ni_1}F_{i_1i_2}\cdots F_{i_s,-n}\ts
\frac{(f_{-n}-f_{j_1}-1)\cdots (f_{-n}-f_{j_k}-1)}{2f_{-n}-2}.
\non
\end{equation}
This formula together with (\ref{zaieq}) is used in the derivation 
of the following relations from
(\ref{zia}) (cf. \cite[Proposition~2.1]{m:br}): for $a=-n,n$
\begin{equation}\label{Fn-1a}
F_{n-1,a}=\sum_{i=-n+1}^{n-1}z_{n-1,i}\ts z_{ia}
\prod_{j=-n+1,\ts j\ne \pm i}^{n-1}\frac{1}{f_i-f_j},
\end{equation}
where $z_{n-1,n-1}:=1$ and the equalities are considered 
in $\U'(\g_n)$ modulo
the ideal $\J$.

Let us now introduce the $\gl(2)$-{\it Yangian\/} $\Y(2)$ and the 
({\it orthogonal\/})
{\it twisted Yangian\/} $\Y^+(2)$; see \cite{mno:yc} for more details.
The Yangian $\Y(2)$ is the
complex associative algebra with the
generators $t_{ab}^{(1)},t_{ab}^{(2)},\dots$ where 
$a,b\in\{-n,n\}$,
and the defining relations
\begin{equation}\label{rel}
[t_{ab}(u),t_{cd}(v)]=\frac{1}{u-v}
\Big(t_{cb}(u)t_{ad}(v)-t_{cb}(v)t_{ad}(u)\Big),
\end{equation}
where
\begin{equation}%\label{ser}
t_{ab}(u): = \delta_{ab} + t^{(1)}_{ab} u^{-1} + t^{(2)}_{ab}u^{-2} +
\cdots \in \Y(2)[[u^{-1}]].
\non
\end{equation}
Introduce the series $s_{ab}(u)$, $a,b\in\{-n,n\}$ by
\begin{equation}%\label{sab}
s_{ab}(u)=t_{an}(u)t_{-b,-n}(-u)+t_{a,-n}(u)t_{-b,n}(-u).
\non
\end{equation}
Write
$
s_{ab}(u)=\delta_{ab}+s_{ab}^{(1)}u^{-1}+s_{ab}^{(2)}u^{-2}+\cdots.
$
The twisted Yangian $\Y^+(2)$ is defined as the subalgebra of $\Y(2)$
generated by the elements $s_{ab}^{(1)},s_{ab}^{(2)},\dots$ where 
$a,b\in\{-n,n\}$. Note that $\Y^+(2)$ can be equivalently defined as
an abstract algebra with these generators and certain linear and quadratic
defining relations; see \cite[Section~3]{mno:yc}.

The Yangian $\Y(2)$ is a Hopf algebra with the 
coproduct 
\begin{equation}\label{cop}
\Delta (t_{ab}(u))=t_{an}(u)\ot
t_{nb}(u)+t_{a,-n}(u)\ot t_{-n,b}(u).
\end{equation}
The twisted Yangian $\Y^+(2)$ is a left coideal in $\Y(2)$
with
\begin{equation}\label{cops}
\Delta (s_{ab}(u))=\sum_{c,d\in\{-n,n\}}t_{ac}(u)t_{-b,-d}(-u)\ot
s_{cd}(u).
\end{equation}

Given a pair of complex numbers $(\alpha,\beta)$ 
such that $\alpha-\beta\in\ZZ_+$
we denote by 
$L(\alpha,\beta)$ the irreducible representation of the Lie algebra
$\gl(2)$ with the highest weight $(\alpha,\beta)$ with respect to the
upper triangular Borel subalgebra. We have 
$\dim L(\alpha,\beta)=\alpha-\beta+1$. We may regard $L(\alpha,\beta)$ as
a $\Y(2)$-module by using the algebra homomorphism $\Y(2)\to\U(\gl(2))$
given by
\begin{equation}\label{hom}
t_{ab}(u)\mapsto \delta_{ab}+E_{ab}u^{-1},\qquad a,b\in\{-n,n\}.
\end{equation}
The coproduct (\ref{cop}) allows one to construct representations of
$\Y(2)$ of the form
\begin{equation}%\label{tenpr}
L=L(\alpha_1,\beta_1)\ot\cdots\ot L(\alpha_k,\beta_k).
\non
\end{equation}
For any $\gamma\in\C$ denote by $W(\gamma)$ the one-dimensional
representation of $\Y^+(2)$ spanned by a vector $w$ such that
\begin{equation}%\label{onedim}
s_{nn}(u)\ts w=\frac{u+\gamma}{u+1/2}\ts w,\qquad 
s_{-n,-n}(u)\ts w=\frac{u-\gamma+1}{u+1/2}\ts w,
\non
\end{equation}
and $s_{a,-a}(u)\ts w=0$ for $a=-n,n$. 
By (\ref{cops}) we can regard the tensor
product $L\ot W(\gamma)$ as a representation of $\Y^+(2)$.
Representations of this type
essentially exhaust all finite-dimensional irreducible representations of
$\Y^+(2)$ \cite{m:fd}. The vector space isomorphism
\begin{equation}\label{isomL}
L\ot W(\gamma)\to L,\qquad v\ot w\mapsto v,\qquad v\in L
\end{equation}
provides $L\ot W(\gamma)$ with an action of $\Y(2)$.

\section{Construction of the basis}
\setcounter{equation}{0}

Introduce the following series with coefficients in the
transvector algebra $\Z(\g_n,\g_{n-1})$: for $a,b\in\{-n,n\}$
\begin{multline}\label{Zab}
Z_{ab}(u)=\Bigg(-\Big(\delta_{ab}(u-n+3/2)+F_{ab}\Big)
\prod_{i=-n+1}^{n-1}(u+g_i)\\
+\sum_{i=-n+1}^{n-1}z_{ai}z_{ib}\ts (u+g_{-i})
\prod_{j=-n+1,\ts j\ne \pm i}^{n-1}\frac{u+g_j}{g_i-g_j}\Bigg)
\frac{1}{2u+1},
\end{multline}
where $g_i:=f_i+1/2$ for all $i$.

As we shall see below (Corollary~\ref{cor:bra})
the space $V(\lambda)^+_{\mu}$ is nonzero
only if
there exist $\nu_1,\dots,\nu_{n-1}$ such that
the inequalities (\ref{ineq}) hold.
We shall be assuming that this condition is satisfied.

\bpr\label{prop:map}
{\rm (i)} The mapping
\begin{equation}\label{map}
s_{ab}(u)\mapsto -2u^{-2n+2}\ts Z_{ab}(u),\qquad a,b\in\{-n,n\}
\end{equation}
defines an algebra homomorphism $\Y^+(2)\to \Z(\g_n,\g_{n-1})$.

{\rm (ii)} The representation of $\Y^+(2)$ in the space
$V(\lambda)^+_{\mu}$ defined via the homomorphism (\ref{map})
is irreducible.
\epr

\Proof We use the same arguments as for the proof of the corresponding
statements in the symplectic case; see \cite[Section~5]{m:br}.
So we shall only give a few key formulas; the details
can be restored by using \cite{m:br}.

Introduce the
$2n\times 2n$-matrix $F=(F_{ij})$ whose $ij$th entry is
the element $F_{ij}\in\g_n$ and set
\begin{equation}%\label{Fu}
 F(u)=1+\frac{F}{u+1/2}.
\non
\end{equation}
Denote by $\wh { F}(u)$ the corresponding {\it Sklyanin comatrix\/};
see \cite[Section~2]{m:fd}. The mapping
\begin{equation}\label{homom}
s_{ab}(u)\mapsto  c(u)\ts  \wh { F}(-u+n-1)_{ab},\qquad a,b\in\{-n,n\},
\end{equation}
where
\begin{equation}%\label{cu}
c(u)=\prod_{k=1}^{n-1} (1-(k-1/2)^2\ts u^{-2}),
\non
\end{equation}
defines an algebra homomorphism from $\Y^+(2)$ to the centralizer
$\CC_n$ of $\g_{n-1}$ in $\U(\g_n)$ \cite[Proposition~2.1]{m:fd};
cf. \cite{o:ty}.
Further, a slight
generalization of \cite[Proposition~3.1]{m:br}
implies the following expression for the
$ab$-entries of the matrix $\wh{ F}(u+1/2)$:
\begin{equation}%\label{quasi}
\wh{ F}(u+1/2)_{ab}=\left(
\delta_{ab}-\sum_{k=1}^{\infty}F^{(k)}_{ab}\ts u^{-k}\right)
\cdot\sdet  F^{(n-1)}(u-1/2).
\non
\end{equation}
Here $\sdet  F^{(n-1)}(u)$ is the
{\it Sklyanin determinant\/} of 
the matrix obtained from $F(u)$ by deleting the 
$(\pm n)$th rows and columns
\cite[Section~2]{m:fd},
and
\begin{equation}%\label{nota}
F^{(k)}_{ab}=\sum F^{}_{ai_1}F^{}_{i_1i_2}\cdots F^{}_{i_{k-1}b},
\non
\end{equation}
summed over the indices 
$i_m\in\{-n+1,\dots,n-1\}$.
Finally, calculating the images of $F^{(k)}_{ab}$ and $\sdet  F^{(n-1)}(u)$
with respect to the natural homomorphism $\pi:\CC_n\to\Z(\g_n,\g_{n-1})$
(cf. \cite[Section~5]{m:br}), we find that the composition
of $\pi$ and (\ref{homom}) yields (\ref{map}).
\endproof

The next theorem provides an identification of the
$\Y^+(2)$-module $V(\lambda)^+_{\mu}$.

\bth\label{thm:isom}
We have an isomorphism of
$\ts\Y^+(2)$-modules
\begin{equation}\label{isom}
V(\lambda)^+_{\mu}\simeq
L(\alpha_1,\beta_1)\ot\cdots\ot L(\alpha_{n-1},\beta_{n-1})
\ot W(-\alpha_0)
\end{equation}
where $\alpha_1=\min\{-|\lambda_{1}|,-|\mu_{1}|\}-1/2$,\quad
$\alpha_0=\alpha_1+|\lambda_1+\mu_1|$,
\begin{align}%\label{alphai}
\alpha_i&=\min\{\lambda_{i},\mu_{i}\}-i+1/2,&\qquad i&=2,\dots,n-1,
\non
\\
\beta_i&=\max\{\lambda_{i+1},\mu_{i+1}\}-i+1/2,&\qquad i&=1,\dots,n-1
\non
\end{align}
with $\mu_{n}:=-\infty$.
In particular, $V(\lambda)^+_{\mu}$ is equipped with
an action of $\Y(2)$ defined by (\ref{isomL}).
\eth

\Proof Consider the following vector in $V(\lambda)^+_{\mu}$
\begin{equation}\label{hv}
\xi_{\mu}=\prod_{i=1}^{n-1}\Big(z_{ni}^{\max\{\lambda_i,\mu_i\}-\mu_i}
z_{i,-n}^{\max\{\lambda_i,\mu_i\}-\lambda_i}\Big)\ts \xi.
\end{equation}
Repeating the arguments of the proof of Theorem~5.2 in \cite{m:br}
we show that $\xi_{\mu}$ is the highest vector of the 
$\Y^+(2)$-module $V(\lambda)^+_{\mu}$. That is, $\xi_{\mu}$
is annihilated by $s_{-n,n}(u)$, and $\xi_{\mu}$ is an eigenvector
for $s_{nn}(u)$. Namely, $s_{nn}(u)\xi_{\mu}=\mu(u)\xi_{\mu}$,
where the highest weight $\mu(u)$ is given by
\begin{equation}\label{muu}
\mu(u)=(1-\alpha_0u^{-1})\cdots(1-\alpha_{n-1}u^{-1})
(1+\beta_1u^{-1})\cdots
(1+\beta_{n-1}u^{-1})(1+\frac12u^{-1})^{-1}.
\end{equation}
This is proved simultaneously with the following relations
by induction on the degree of the monomial in (\ref{hv}):
for $i=1,\dots,n-1$
\begin{multline}%\label{zinact}
z_{in}\ts\xi_{\mu}=-(m_i+\wt{\alpha}_1+1)\cdots \wh{(m_i+\alpha_i+1)}
\cdots(m_i+\alpha_{n-1}+1)\\
{}\times (m_i-\beta_0+1)\cdots (m_i-\beta_{n-1}+1)\ts \xi_{\mu+\delta_i},
\non
\end{multline}
and
\begin{equation}%\label{z-niact}
z_{-ni}\ts\xi_{\mu}=-(m_i-\wt{\alpha}_1)\cdots(m_i-\alpha_{n-1})
(m_i+\beta_0)\cdots \wh{(m_i+\beta_{i-1})}\cdots 
(m_i+\beta_{n-1})\ts \xi_{\mu-\delta_i},
\non
\end{equation}
where we have used the notation 
\begin{align}%\label{notmal}
m_i&=\mu_i-i+1/2,\qquad i=1,\dots,n-1,\non\\
\wt{\alpha}_1&=\min\{\lambda_1,\mu_1\}-1/2,\qquad
\beta_0=\max\{\lambda_1,\mu_1\}+1/2,\non
\end{align}
(note that $\{\wt{\alpha}_1,-\beta_0\}=\{\alpha_0,\alpha_1\}$).
On the other hand, it follows from \cite[Corollary~6.6]{m:fd} that the
tensor product in (\ref{isom}) is an irreducible representation
of $\Y^+(2)$. Its highest weight can be easily calculated and is given by
the same formula (\ref{muu}). \endproof

Set
$
T_{ab}(u)=u^{n-1}\ts t_{ab}(u)
$
for $a,b\in\{-n,n\}$.
By (\ref{cop}) and (\ref{hom}), $T_{ab}(u)$, 
as an operator in $V(\lambda)^+_{\mu}$,
is a polynomial in $u$:
\begin{equation}\label{Tabpol}
T_{ab}(u)=\delta_{ab}u^{n-1}+t_{ab}^{(1)}u^{n-2}+\cdots+t_{ab}^{(n-1)}.
\end{equation}
By (\ref{rel}), (\ref{cops}) and (\ref{map})
we have an equality of operators in $V(\lambda)^+_{\mu}$:
\begin{equation}\label{ZT}
Z_{n,-n}(u)=\frac{(u-\alpha_0)T_{n,-n}(-u)T_{nn}(u)+
(u+\alpha_0)T_{n,-n}(u)T_{nn}(-u)}{(-1)^n\ts 2u}.
\end{equation}
Therefore, $Z_{n,-n}(u)$ is a polynomial in $u^2$ of degree $n-2$.
On the other hand, we find from (\ref{Zab}) that 
$Z_{n,-n}(-g_i)=z_{ni}z_{i,-n}$. Thus, 
by the Lagrange interpolation formula,
$Z_{n,-n}(u)$ can also be given by
\begin{equation}\label{Zn-n}
Z_{n,-n}(u)=\sum_{i=1}^{n-1}z_{ni}z_{i,-n}\prod_{j=1,\ts j\ne i}^{n-1}
\frac{u^2-g_j^2}{g_i^2-g_j^2}.
\end{equation}

\noindent
{\it Remark.} To make the above evaluation $Z_{n,-n}(-g_i)$
well-defined we agree to consider the series $Z_{ab}(u)$
with $a,b\in\{-n,n\}$
as elements of the {\it right\/} module over the field
of rational functions in $\h$ and $u$ generated by
monomials in the $z_{ia}$. \endproof

Theorem~3.2 implies that basis vectors of $V(\lambda)^+_{\mu}$ can be
naturally parametrized by $(n-1)$-tuples $(\nu_1,\dots,\nu_{n-1})$, where
all the entries are simultaneously integers or half-integers
together with the $\lambda_i$ and the $\mu_i$,
and the following inequalities hold:
\begin{equation}\label{ineq}
\begin{split}
-|\lambda_1|&\geq\nu_1\geq\lambda_2\geq\nu_2\geq\lambda_3\geq
\cdots \geq\lambda_{n-1}\geq\nu_{n-1}\geq\lambda_n,\\
-|\mu_1|&\geq\nu_1\geq\mu_2\geq\nu_2\geq\mu_3\geq
\cdots \geq\mu_{n-1}\geq\nu_{n-1}.
\end{split}
\end{equation}
For $i\geq 1$ set
\begin{equation}%\label{gammal}
\gamma_i=\nu_i-i+1/2,\qquad\qquad l_i=\lambda_{i}-i+1/2.
\non
\end{equation}
Introduce the vectors
\begin{equation}%\label{xinumu}
\xi_{\nu\mu}=\prod_{i=1}^{n-1}Z_{n,-n}(\gamma_i-1)\cdots 
Z_{n,-n}(\beta_i+1)Z_{n,-n}(\beta_i)\ts\xi_{\mu}.
\non
\end{equation}
Using (\ref{Zn-n}) we can write an equivalent expression; 
cf. \cite[Section~6]{m:br}:
\begin{equation}\label{xinumuz}
\xi_{\nu\mu}=\prod_{i=1}^{n-1}z_{ni}^{\nu_{i-1}-\mu_i}
z_{i,-n}^{\nu_{i-1}-\lambda_i}\cdot
\prod_{k=l_{n}+1}^{\gamma_{n-1}-1}Z_{n,-n}(k)\ts\xi,
\end{equation}
where $\nu_0:=\max\{\lambda_1,\mu_1\}$. 
The vectors $\xi_{\nu\mu}$ with $\nu$ satisfying (\ref{ineq})
form a basis of the space
$V(\lambda)^+_{\mu}$;
see \cite[Proposition~6.1]{m:br}.
We shall use the following
normalized basis vectors
\begin{equation}%\label{zetanumu}
\zeta_{\nu\mu}=\prod_{1\leq i<j\leq n-1}
(-\gamma_i-\gamma_j)!\  \xi_{\nu\mu}.
\non
\end{equation}
The generators of the Yangian $\Y(2)$
act in the basis $\{\zeta_{\nu\mu}\}$
by the rule: for $i=1,\dots,n-1$
\begin{align}
T_{nn}(u)\ts\zeta_{\nu\mu}&=(u+\gamma_1)\cdots (u+\gamma_{n-1})
\ts\zeta_{\nu\mu}, \non\\
T_{n,-n}(-\gamma_i)\ts \zeta_{\nu\mu}&=
\frac{1}{\gamma_i-\alpha_0}\ts
\zeta_{\nu+\delta_i,\mu},\label{acT1}
\\
T_{-n,n}(-\gamma_i)\ts \zeta_{\nu\mu}&=
\prod_{k=0}^{n-1}(\alpha_k-\gamma_i+1)\prod_{k=1}^{n-1}(\beta_k-\gamma_i)
\ts\zeta_{\nu-\delta_i,\mu};\non
\end{align}
cf. \cite[Proposition~4.2]{m:br}.
The action of $T_{-n,-n}(u)$ can be found by using the {\it
quantum determinant\/} 
\begin{align}\label{qdet1}
d(u)&=T_{-n,-n}(u+1)T_{nn}(u)-T_{n,-n}(u+1)T_{-n,n}(u)\\
\label{qdet2}
{}&=T_{-n,-n}(u)T_{nn}(u+1)-T_{-n,n}(u)T_{n,-n}(u+1);
\end{align}
see, e.g. \cite[Section~2]{mno:yc}.
The coefficients of the quantum determinant belong to the center of $\Y(2)$
and so, $d(u)$ acts in $V(\lambda)^+_{\mu}$ as a scalar which
can be found by the application of (\ref{qdet1}) to the 
highest weight vector $\xi_{\mu}$. 
So, we have
\begin{equation}%\label{qdetact}
d(u)\ts \zeta_{\nu\mu}=(u+\alpha_1+1)\cdots (u+\alpha_{n-1}+1)
(u+\beta_1)\cdots(u+\beta_{n-1})\ts \zeta_{\nu\mu}.
\non
\end{equation}
Now, using (\ref{acT1}) and (\ref{qdet2}) we obtain
\begin{multline}\label{acT-n}
T_{-n,-n}(u)\ts \zeta_{\nu\mu}=\prod_{i=1}^{n-1}
\frac{(u+\alpha_i+1)(u+\beta_i)}{u+\gamma_i+1}\ts \zeta_{\nu\mu}\\
+\prod_{i=1}^{n-1}
\frac{1}{u+\gamma_i+1}\ts T_{-n,n}(u)T_{n,-n}(u+1)\ts \zeta_{\nu\mu}.
\end{multline}
The operators $T_{-n,n}(u)$ and $T_{n,-n}(u)$ are polynomials in $u$
of degree ${}\leq n-2$; see (\ref{Tabpol}). 
Therefore, their action can be found from
(\ref{acT1}) by using the Lagrange interpolation formula.

The following branching 
rule for the reduction $\g_n\downarrow\g_{n-1}$ is implied by
Theorem~\ref{thm:isom}; cf. \cite[Corollary~5.3]{m:br}.

\bco\label{cor:bra}
The restriction of $V(\lambda)$ to the subalgebra $\g_{n-1}$
is isomorphic to the direct sum
$\bigoplus\ts c(\mu)V'(\mu)$
of finite-dimensional irreducible 
representations $V'(\mu)$ of $\g_{n-1}$
where the multiplicity $c(\mu)$ equals the number of $(n-1)$-tuples
$\nu$ satisfying the inequalities (\ref{ineq}).
\eco

\Proof We have $c(\mu)=\dim V(\lambda)^+_{\mu}$. By Theorem~\ref{thm:isom},
\begin{equation}%\label{dim}
\dim V(\lambda)^+_{\mu}=\prod_{i=1}^{n-1}(\alpha_i-\beta_i+1),
\non
\end{equation}
if there exists $\nu$ satisfying (\ref{ineq}). Otherwise,
the space $V(\lambda)^+_{\mu}$ is trivial. This is proved
by comparison of the dimensions of $V(\lambda)$ and 
$\bigoplus\ts c(\mu)V'(\mu)$ 
with the use of
\cite[Chapter VII, Section 9]{w:cg}. 
\endproof

Applying the above construction of the vectors $\zeta_{\nu\mu}$ to
the subalgebras of the chain
\begin{equation}%\label{chaing}
\g_1\subset\g_2\subset\cdots\subset\g_n,\qquad \g_k=\oa(2k)
\non
\end{equation}
we obtain a basis of $V(\lambda)$ parametrized by the
$D$-{\it type Gelfand--Tsetlin patterns\/}
(cf.~\cite{l:cc}) which we denote by $\Lambda$:
\begin{align}
&\lambda^{}_{n1}\qquad\lambda^{}_{n2}
\qquad\qquad\qquad\cdots\qquad\qquad\qquad\lambda^{}_{nn}\non\\
&\qquad \lambda'_{n-1,1}
\qquad\qquad\cdots\qquad\qquad\lambda'_{n-1,n-1}\non\\
&\lambda^{}_{n-1,1}\qquad\cdots\qquad\lambda^{}_{n-1,n-1}\non\\
&\qquad\qquad\cdots\qquad\cdots\non\\
&\lambda^{}_{21}\qquad\lambda^{}_{22}\non\\
&\qquad\lambda'_{11}\non\\
&\lambda^{}_{11}\non
\end{align}
Here 
the upper row coincides with $\lambda$, 
all the entries are simultaneously integers or half-integers
and 
the following inequalities hold
\begin{align}
-|\lambda^{}_{k1}|&\geq\lambda'_{k-1,1}\geq\lambda^{}_{k2}\geq
\lambda'_{k-1,2}\geq \cdots\geq
\lambda^{}_{k,k-1}\geq\lambda'_{k-1,k-1}\geq\lambda^{}_{kk},\non\\
-|\lambda^{}_{k-1,1}|&\geq\lambda'_{k-1,1}\geq\lambda^{}_{k-1,2}\geq
\lambda'_{k-1,2}\geq \cdots\geq
\lambda^{}_{k-1,k-1}\geq\lambda'_{k-1,k-1}\non
\end{align}
for $k=2,\dots,n$. Set
\begin{equation}\label{ll}
l^{}_{ki}=\lambda^{}_{ki}-i+1,\qquad l'_{ki}=\lambda'_{ki}-i+1,\qquad
1\leq i\leq k\leq n
\end{equation}
and introduce the vectors
\begin{equation}%\label{basisxi}
\xi^{}_{\Lambda}=
\prod_{k=2,\dots,n}^{\rightarrow}
\left(\prod_{i=1}^{k-1}z_{ki}^{\lambda'_{k-1,i-1}-\lambda^{}_{k-1,i}}
z_{i,-k}^{\lambda'_{k-1,i-1}-\lambda^{}_{ki}}
\prod_{q=l^{}_{kk}+1}^{l'_{k-1,k-1}-1}
Z_{k,-k}(q-\frac12)\right)\xi
\non
\end{equation}
with $\lambda'_{k-1,0}:=\max\{\lambda^{}_{k1},\lambda^{}_{k-1,1}\}$. 
Finally, set
\begin{equation}%\label{basiszeta}
\zeta^{}_{\Lambda}=N^{}_{\Lambda}\ts \xi^{}_{\Lambda},\qquad
N^{}_{\Lambda}=\prod_{k=2}^{n-1}\prod_{1\leq i<j\leq k}(-l'_{ki}-l'_{kj}+1)!
\non
\end{equation}
The following proposition is implied by
Corollary~\ref{cor:bra}.

\bpr\label{prop:basis}
The vectors $\zeta^{}_{\Lambda}$ parametrized by the Gelfand--Tsetlin 
patterns $\Lambda$ form a basis of the representation $V(\lambda)$. \endproof
\epr

\section{Matrix element formulas}
\setcounter{equation}{0}

Introduce the
following elements of $\U(\g_n)$:
\begin{equation}%\label{Phi}
\Phi_{-k,k}=\sum_{i=1}^{k-1}F_{-k,i}F_{ik},\qquad k=2,\dots,n.
\non
\end{equation}
We shall find the action of $\Phi_{-k,k}$ in the basis $\{\zeta^{}_{\Lambda}\}$,
which will be used later on.
Since $\Phi_{-k,k}$ commutes with the subalgebra $\g_{k-1}$ it suffices
to consider the case $k=n$. The image of $\Phi_{-n,n}$ under
the natural homomorphism $\pi:\CC_n\to\Z(\g_n,\g_{n-1})$
coincides with the coefficient at $u^{2n-4}$ of the polynomial
$Z_{-n,n}(u)$; see the proof of Proposition~\ref{prop:map}.
The following analog of (\ref{ZT}) is obtained from
(\ref{rel}), (\ref{cops}) and (\ref{map}):
\begin{equation}%\label{ZT2}
Z_{-n,n}(u)=
\frac{(u-\alpha_0)T_{-n,-n}(-u)T_{-n,n}(u)+
(u+\alpha_0)T_{-n,-n}(u)T_{-n,n}(-u)}{(-1)^n\ts 2u}.
\non
\end{equation}
Therefore, we have an equality
of operators in $V(\lambda)^+_{\mu}$:
\begin{equation}\label{Phit}
\Phi_{-n,n}=-t^{(2)}_{-n,n}+t^{(1)}_{-n,n}t^{(1)}_{-n,-n}+
(1+\alpha_0)\ts t^{(1)}_{-n,n}.
\end{equation}
The image of $s_{nn}^{(1)}$ under the homomorphism (\ref{homom})
is $F_{nn}$. On the other hand, by (\ref{cops}) we have
\begin{equation}%\label{snn1}
s_{nn}^{(1)}=t_{nn}^{(1)}-t_{-n,-n}^{(1)}-\alpha_0-1/2,
\non
\end{equation}
as operators in $V(\lambda)^+_{\mu}$. Therefore, (\ref{Phit}) can be written as
\begin{equation}%\label{Phit2}
\Phi_{-n,n}=-t^{(2)}_{-n,n}+t^{(1)}_{-n,n}t^{(1)}_{nn}-
(F_{nn}+3/2)\ts t^{(1)}_{-n,n}.
\non
\end{equation}
Finally, relations (\ref{acT1}) imply that
\begin{equation}\label{Phiact}
\Phi_{-n,n}\zeta_{\nu\mu}=\sum_{i=1}^{n-1}\theta_i\ts(F_{nn}-\gamma_i+3/2)
\ts \zeta_{\nu-\delta_i,\mu},
\end{equation}
where
\begin{equation}%\label{theta}
\theta_i=-\prod_{k=0}^{n-1}(\alpha_k-\gamma_i+1)
\prod_{k=1}^{n-1}(\beta_k-\gamma_i)
\prod_{j=1,\ts j\ne i}^{n-1}(\gamma_j-\gamma_i)^{-1}.
\non
\end{equation}

The action of $F_{nn}$ in $V(\lambda)^+_{\mu}$ is immediately found so that
\begin{equation}%\label{Fnn}
F_{nn}\ts\zeta_{\nu\mu}=\left(2\ts\sum_{i=0}^{n-1}\nu_i-
\sum_{i=1}^n\lambda_i-\sum_{i=1}^{n-1}\mu_i\right)\zeta_{\nu\mu}.
\non
\end{equation}

\noindent
{\it Remark.} One can introduce the elements $\Phi_{k,-k}$ by
\begin{equation}\non
\Phi_{k,-k}=\sum_{i=1}^{k-1}F_{ki}F_{i,-k}.
\end{equation}
The action of $\Phi_{n,-n}$ on $V(\lambda)^+_{\mu}$ is found in the same
way as that of $\Phi_{-n,n}$:
\begin{equation}\non
\Phi_{n,-n}\zeta_{\nu\mu}=\sum_{i=1}^{n-1}
\prod_{j=1,\ts j\ne i}^{n-1}\frac{1}{\gamma_j-\gamma_i}
\ts \zeta_{\nu+\delta_i,\mu},
\end{equation}
although this will not be used. \endproof

The operator
$F_{n-1,-n}$ preserves
the subspace of $\g_{n-2}$-highest vectors in $V(\lambda)$.
Therefore it suffices to calculate its action
on the basis vectors of the form
\begin{equation}\label{numunu}
\xi_{\nu\mu\nu'}=X_{\mu\nu'}\ts\xi_{\nu\mu},
\end{equation}
where 
$X_{\mu\nu'}$ denotes the operator
\begin{equation}%\label{X}
X_{\mu\nu'}=\prod_{i=1}^{n-2}
z_{n-1,i}^{\nu'_{i-1}-\mu'_i}\ts z_{i,-n+1}^{\nu'_{i-1}-\mu_i}\cdot
\prod_{a=m_{n-1}+1}^{\gamma'_{n-1}-1}
Z_{n-1,-n+1}(a),
\non
\end{equation}
$\nu'$ and $\mu'$ are $(n-2)$-tuples of integers or half-integers
such that the inequalities (\ref{ineq})
are satisfied with $\lambda$, $\nu$, $\mu$ respectively replaced
by $\mu$, $\nu'$, $\mu'$; we set $\gamma'_i=\nu'_i-i+1/2$
and $\nu'_0=\max\{\mu_1,\mu'_1\}$. 
The operator $F_{n-1,-n}$ is permutable with the 
elements $z_{n-1,i}$, $z_{i,-n+1}$ and
$Z_{n-1,-n+1}(u)$ which follows from their explicit formulas. 
Hence, we can write
\begin{equation}\label{FX}
F_{n-1,-n}\ts \xi_{\nu\mu\nu'}=X_{\mu\nu'}\ts F_{n-1,-n}\ts \xi_{\nu\mu}.
\end{equation}
By (\ref{Fn-1a}) (with $a=-n$) we need to express
\begin{equation}\label{Xgen}
X_{\mu\nu'}\ts z_{n-1,i}z_{i,-n}\ts \xi_{\nu\mu},\qquad
i=-n+1,\dots,n-1
\end{equation}
as a linear combination of the vectors $\xi_{\nu\mu\nu'}$.
If $i\ne\pm 1$ then the calculation is exactly the same as in
\cite[Section~6]{m:br} where one uses the relations
\begin{equation}\label{eval}
Z_{n,-n}(-g_n)=z_{n,-n},\qquad 
Z_{n,-n}(-g_i)=z_{ni}z_{i,-n},
\end{equation}
which follow from (\ref{zn-n}) and (\ref{Zn-n}).
Now consider (\ref{Xgen}) with $i=-1$. We have
\begin{equation}%\label{XZ-1}
X_{\mu\nu'}\ts z_{n-1,-1}z_{-1,-n}\ts \xi_{\nu\mu}=
-X_{\mu\nu'}\ts z_{1,-n+1}z_{n1}\ts \xi_{\nu\mu}.
\non
\end{equation}
If $\lambda_1\geq \mu_1$ then 
$z_{n1}\ts \xi_{\nu\mu}=\xi_{\nu,\mu-\delta_1}$
while for $\lambda_1< \mu_1$ we derive from (\ref{eval})
that
\begin{equation}%\label{z-1xile}
z_{n1}\ts \xi_{\nu\mu}=\sum_{i=1}^{n-1}\prod_{a=1,\ts a\ne i}^{n-1}
\frac{m_1^2-\gamma_a^2}
{\gamma_i^2-\gamma_a^2}\ts\xi_{\nu+\delta_i,\mu-\delta_1}.
\non
\end{equation}
Similarly, if $\mu'_1\geq \mu_1$ then
$
X_{\mu\nu'}\ts z_{1,-n+1}=X_{\mu-\delta_1,\nu'}
$
while for $\mu'_1< \mu_1$ one has
\begin{equation}%\label{Xzpartle}
X_{\mu\nu'}\ts z_{1,-n+1}=\sum_{r=1}^{n-2}\prod_{a=1,\ts a\ne r}^{n-2}
\frac{m_1^2-{\gamma'_a}^2}{{\gamma'_r}^2-{\gamma'_a}^2}\ts
X_{\mu-\delta_1,\nu'+\delta_r}.
\non
\end{equation}
Finally, take $i=1$ in (\ref{Xgen}). If $\lambda_1\leq \mu_1$ then
$
z_{1,-n}\ts \xi_{\nu\mu}=\xi_{\nu,\mu+\delta_1},
$
and if $\lambda_1> \mu_1$ then
\begin{equation}%\label{z1xige}
z_{1,-n}\ts \xi_{\nu\mu}=
\sum_{i=1}^{n-1}\prod_{a=1,\ts a\ne i}^{n-1}\frac{(m_1+1)^2-\gamma_a^2}
{\gamma_i^2-\gamma_a^2}\ts \xi_{\nu+\delta_i,\mu+\delta_1}.
\non
\end{equation}
Similarly, if $\mu'_1\leq \mu_1$ then
$
X_{\mu\nu'}\ts z_{n-1,1}=X_{\mu+\delta_1,\nu'}
$
and if $\mu'_1> \mu_1$ then
\begin{equation}%\label{Xzpartlen}
X_{\mu\nu'}\ts z_{n-1,1}=\sum_{r=1}^{n-2}\prod_{a=1,\ts a\ne r}^{n-2}
\frac{(m_1+1)^2-{\gamma'_a}^2}{{\gamma'_r}^2-{\gamma'_a}^2}\ts
X_{\mu+\delta_1,\nu'+\delta_r}.
\non
\end{equation}
 
The action of the elements $F_{n-1,n}$ on the vectors (\ref{numunu})
can be expressed in two different ways. First we sketch
a calculation similar to the one used above which leads to
(rather complicated) explicit formulas for the matrix elements. 
Then we give slightly less explicit but more convenient
formulas where $F_{n-1,n}$ is represented by a commutator-like expression
of simpler operators.

We have the following analog of (\ref{FX}):
\begin{equation}\label{FXn}
F_{n-1,n}\ts \xi_{\nu\mu\nu'}=X_{\mu\nu'}\ts F_{n-1,n}\ts \xi_{\nu\mu}.
\end{equation}
Now use (\ref{Fn-1a}) with $a=n$. Here we need to calculate
$z_{in}\ts \xi_{\nu\mu}$ instead of $z_{i,-n}\ts \xi_{\nu\mu}$
in the previous case. Suppose that $i>1$. We have
\begin{equation}%\label{zinxi}
z_{in}\ts \xi_{\nu\mu}=z_{in}z_{ni}\ts\xi_{\nu,\mu+\delta_i}=
z_{-n,-i}z_{-i,-n}\ts\xi_{\nu,\mu+\delta_i}=
Z_{-n,-n}(-g_{-i})\ts\xi_{\nu,\mu+\delta_i};
\non
\end{equation}
see (\ref{Zab}). However, 
\begin{equation}%\label{gixi}
-g_{-i}\ts \xi_{\nu,\mu+\delta_i}=(m_i+1)\ts \xi_{\nu,\mu+\delta_i}.
\non
\end{equation}
To calculate $Z_{-n,-n}(m_i+1)\ts\xi_{\nu,\mu+\delta_i}$ we use the
following equality of operators in $V(\lambda)^+_{\mu}$:
\begin{equation}%\label{ZT22}
Z_{-n,-n}(u)=
\frac{(u-\alpha_0)T_{-n,n}(u)T_{n,-n}(-u)+
(u+\alpha_0+1)T_{-n,-n}(u)T_{nn}(-u)}{(-1)^n\ts (2u+1)},
\non
\end{equation}
see (\ref{rel}), (\ref{cops}) and (\ref{map}); and then apply
formulas (\ref{acT1}) and (\ref{acT-n}).

To calculate $z_{-i,n}\ts \xi_{\nu\mu}$
we first permute $z_{-i,n}$ with the generators $z_{nj}$ and $z_{j,-n}$ with
$j=1,\dots,i-1$
in (\ref{xinumuz}). Further, we use the relation
\begin{equation}%\label{zzZ}
z_{-i,n}z_{i,-n}=(-1)^{n-i}z_{-n,i}z_{i,-n}=(-1)^{n-i}Z_{-n,-n}(-g_{i})
\non
\end{equation}
and complete the calculation in a similar manner.
To find $z_{\pm 1,n}\ts \xi_{\nu\mu}$ we need to consider
a few different cases which depend on the relationship between
the parameters $\lambda_1$, $\mu_1$ and $\mu'_1$ and then proceed
exactly as above in the calculation of the action of $F_{n-1,-n}$.

We now give an alternative way
of computing the action of 
$F_{n-1,n}$. The basic idea is to replace the operator $z_{in}$
in the above calculation of $z_{in}\ts \xi_{\nu\mu}$ by
the following expression:
for $i=-n+1,\dots,n-1$
\begin{equation}\label{zincom}
z_{in}=[z_{i,-n},\Phi_{-n,n}]\ts \frac{1}{f_i+F_{nn}}
\end{equation}
and then use the formulas for the action of $z_{i,-n}$ and $\Phi_{-n,n}$;
see (\ref{Phiact}). More precisely, we regard (\ref{zincom})
as a relation in the transvector algebra $\Z(\g_n,\g_{n-1})$
which can be proved 
as follows. First, we calculate 
the commutator $[F_{i,-n},\Phi_{-n,n}]$ in $\U(\g_n)$
then consider it modulo the ideal $\J$
and apply the extremal projection $p$
(see Section 2). 

We have
$
\Phi_{-n,n}F_{nn}=(F_{nn}+2)\ts\Phi_{-n,n}
$
and so, (\ref{Fn-1a}), (\ref{FXn}) and (\ref{zincom}) imply that
\begin{equation}\label{Fact}
F_{n-1,n}\ts\xi_{\nu\mu\nu'}
=X_{\mu\nu'}\ts
\left(\Phi_{n-1,-n}(2)\ts\Phi_{-n,n}-
\Phi_{-n,n}\Phi_{n-1,-n}(0)\right)\ts\xi_{\nu\mu},
\end{equation}
where
\begin{equation}%\label{Phinn}
\Phi_{n-1,-n}(u)=\sum_{i=-n+1}^{n-1}z_{n-1,i}\ts z_{i,-n}
\prod_{a=-n+1,\ts a\ne \pm i}^{n-1}\frac{1}{f_i-f_a}\cdot\frac{1}{u+f_i+F_{nn}}.
\non
\end{equation}
The action of $\Phi_{n-1,-n}(u)$
is found exactly as that of $F_{n-1,-n}$.
Note that the operators $X_{\mu\nu'}$ and $\Phi_{-n,n}$ commute.

\noindent
{\it Remark.} The operator $\Phi_{n-1,-n}(u)$ is 
a rational function in $u$ which can have singularities
at the values
$u=0$ and $u=2$ in (\ref{Fact}).
However, the operator
\begin{equation}%\label{operreg}
\Phi_{n-1,-n}(u+2)\ts\Phi_{-n,n}-\Phi_{-n,n}\Phi_{n-1,-n}(u)
\non
\end{equation}
is regular at $u=0$ and coincides with $F_{n-1,n}$. Note the similarity
with the symplectic case \cite{m:br}, where we have 
$2F_{n-1,n}=[F_{n-1,-n}, F_{-n,n}]$. \endproof 

The elements $F^{}_{k-1,-k}$, $F^{}_{k-1,k}$ with $k=2,\dots,n$
and $F^{}_{21}$, $F^{}_{-2,1}$ generate $\g_n$ as a Lie algebra.
Summarizing the above calculations we obtain the following formulas
for the matrix elements of the generators. Given a pattern $\Lambda$
we use the notation (\ref{ll}) and set for $1\leq i<k\leq n$:
\begin{align}
A_{ki}&=\prod_{a=1, \ts a\ne i}^{k-1}
\frac{1}{l_{k-1,i}^2-l_{k-1,a}^2},\non\\
B_{ki}(x)&=\prod_{a=1, \ts a\ne i}^{k-1}
\frac{(x+l'_{k-1,a})(x-l'_{k-1,a}+1)}{l'_{k-1,a}-l'_{k-1,i}},
\non
\end{align}
and
\begin{align}
C_{ki}=&
(\max\{\lambda^{}_{k1},\lambda^{}_{k-1,1}\}+l'_{k-1,i}-1)
(\min\{\lambda^{}_{k1},\lambda^{}_{k-1,1}\}-l'_{k-1,i}+1)\non\\
{}\times&\prod_{a=2}^k (l^{}_{ka}-l'_{k-1,i}+1)
\prod_{a=2}^{k-1} (l^{}_{k-1,a}-l'_{k-1,i}+1)
\prod_{a=1,\ts a\ne i}^{k-1} \frac{1}{l'_{k-1,a}-l'_{k-1,i}}.\non
\end{align}
We denote by
$\Lambda\pm\delta^{}_{ki}$ and $\Lambda\pm\delta'_{ki}$ the arrays
obtained from $\Lambda$ by replacing $\lambda^{}_{ki}$ and $\lambda'_{ki}$
by $\lambda^{}_{ki}\pm1$ and $\lambda'_{ki}\pm1$ respectively. 
Consider the basis $\{\zeta^{}_{\Lambda}\}$ of the representation $V(\lambda)$;
see Proposition \ref{prop:basis}.
We shall suppose
that $\zeta^{}_{\Lambda}=0$ if the array $\Lambda$ is not a pattern.

\bth
The action of the generators of the Lie algebra $\oa(2n)$
in the vasis $\{\zeta^{}_{\Lambda}\}$ is given by the
following formulas.
\begin{align}
F^{}_{kk}\ts \zeta^{}_{\Lambda}=&\left(2\sum_{i=1}^{k}\lambda'_{k-1,i-1}-
\sum_{i=1}^k\lambda^{}_{ki}-\sum_{i=1}^{k-1}\lambda^{}_{k-1,i}\right)
\zeta^{}_{\Lambda},\non\\
F^{}_{k-1,-k}\ts \zeta^{}_{\Lambda}=&\sum_{i=1}^{k-1}A_{ki}
\left(\zeta^+_{\Lambda}(k,i)-\zeta^-_{\Lambda}(k,i)\right).\non
\end{align}
Here
\begin{equation}
\zeta^+_{\Lambda}(k,i)=\sum_{j=1}^{k-1}\sum_{m=1}^{k-2}
B_{kj}(l_{k-1,i})B_{k-1,m}(l_{k-1,i})\ts
\zeta^{}_{\Lambda+\delta'_{k-1,j}+\delta^{}_{k-1,i}+\delta'_{k-2,m}}
\non
\end{equation}
for $i=2,\dots,k-1$; and for $i=1$ if 
$\lambda_{k-1,1}< \lambda_{k1},\lambda_{k-2,1}$.
Otherwise,
\begin{alignat}{2}
\zeta^+_{\Lambda}(k,1)&=\zeta^{}_{\Lambda+\delta^{}_{k-1,1}}
&\qquad&
\text{if} \quad\lambda_{k-1,1}\geq \lambda_{k1},\lambda_{k-2,1},\non\\
{}&=\sum_{j=1}^{k-1}B_{kj}(l_{k-1,1})\ts 
\zeta^{}_{\Lambda+\delta'_{k-1,j}+\delta^{}_{k-1,1}}
&&
\text{if} \quad\lambda_{k-2,1}\leq \lambda_{k-1,1}<\lambda_{k1},\non\\
{}&=\sum_{m=1}^{k-2}B_{k-1,m}(l_{k-1,1})\ts
\zeta^{}_{\Lambda+\delta^{}_{k-1,1}+\delta'_{k-2,m}}
&&
\text{if} \quad\lambda_{k1}\leq \lambda_{k-1,1}<\lambda_{k-2,1}.
\non
\end{alignat}
Furthermore,
\begin{equation}
\zeta^-_{\Lambda}(k,i)=
\zeta^{}_{\Lambda-\delta^{}_{k-1,i}}
\non
\end{equation}
for $i=2,\dots,k-1$; and for $i=1$ if 
$\lambda_{k-1,1}\leq \lambda_{k1},\lambda_{k-2,1}$.
Otherwise,
\begin{alignat}{2}
\zeta^-_{\Lambda}(k,1)&=\sum_{j=1}^{k-1}B_{kj}(l_{k-1,1}-1)\ts 
\zeta^{}_{\Lambda+\delta'_{k-1,j}-\delta^{}_{k-1,1}}
&&
\qquad\quad\text{if} \quad\lambda_{k1}< \lambda_{k-1,1}\leq\lambda_{k-2,1},\non\\
{}&=\sum_{m=1}^{k-2}B_{k-1,m}(l_{k-1,1}-1)\ts
\zeta^{}_{\Lambda-\delta^{}_{k-1,1}+\delta'_{k-2,m}}
&&
\qquad\quad\text{if} \quad\lambda_{k-2,1}< \lambda_{k-1,1}\leq\lambda_{k1},\non\\
{}&=
\sum_{j=1}^{k-1}\sum_{m=1}^{k-2}
B_{kj}(l_{k-1,1}-1)B_{k-1,m}(l_{k-1,1}-1)&&
\zeta^{}_{\Lambda+\delta'_{k-1,j}-\delta^{}_{k-1,1}+\delta'_{k-2,m}}\non\\
&&&
\qquad\quad\text{if} \quad\lambda_{k-1,1}> \lambda_{k1},\lambda_{k-2,1}.
\non
\end{alignat}
The action of $F^{}_{k-1,k}$ is found from the relation
\begin{equation}
F^{}_{k-1,k}=\Bigl[\Phi_{k-1,-k}(u+2)\ts\Phi_{-k,k}-
\Phi_{-k,k}\Phi_{k-1,-k}(u)\Bigr]^{}_{u=0},\non
\end{equation}
where
\begin{equation}
\Phi_{-k,k}\ts\zeta^{}_{\Lambda}=\sum_{i=1}^{k-1} C_{ki}\ts
(F^{}_{kk}-l'_{k-1,i}+2)\ts\zeta^{}_{\Lambda-\delta'_{k-1,i}}\non
\end{equation}
and
\begin{multline}
\Phi_{k-1,-k}(u)\ts \zeta^{}_{\Lambda}=\\
\sum_{i=1}^{k-1}A_{ki}
\Biggl(
\frac{1}{u+l^{}_{k-1,i}+F^{}_{kk}-1}\ts\zeta^+_{\Lambda}(k,i)
-\frac{1}{u-l^{}_{k-1,i}+F^{}_{kk}-1}\ts\zeta^-_{\Lambda}(k,i)
\Biggr).
\non
\end{multline}
\eth
\endproof

\noindent
{\it Example.} Let $n=2$. 
We have $z_{21}=F^{}_{21}$, $z_{1,-2}=F^{}_{1,-2}$ and 
$Z_{2,-2}(u)=F^{}_{21}F^{}_{1,-2}$.
Therefore,
the basis vectors are given by
\begin{equation}
\zeta^{}_{\Lambda}=F_{21}^{\lambda'_{10}-\lambda^{}_{11}}
F_{1,-2}^{\lambda'_{10}-\lambda^{}_{21}}
(F^{}_{21}F^{}_{1,-2})^{\lambda'_{11}-\lambda^{}_{22}}\ts\xi,\non
\end{equation}
where $\lambda'_{10}=\max\{\lambda^{}_{21},\lambda^{}_{11}\}$.
The Lie algebra $\oa(4)$ is isomorphic to the direct sum of two
copies of $\mathfrak{sl}(2)$ and the action of their generators
in the basis $\{\zeta^{}_{\Lambda}\}$
is easily found. 
The resulting formulas also 
hold for the action of the elements of the subalgebra
$\g_2\subset\g_n$ in the basis $\{\zeta^{}_{\Lambda}\}$
of the $\g_n$-module $V(\lambda)$. \endproof

\section*{Acknowledgements}

The author wishes to thank R.~Donnelly, M.~Gould, M.~Moshinsky, and
R.~Proctor for stimulating discussions of the results of \cite{m:br}.

\end{document}